\documentclass[leqno,12pt]{amsart} 
\usepackage{amssymb,amsmath,bbm}
\usepackage{CJKutf8}
\theoremstyle{plain} 
\newtheorem{theorem}{\indent\sc Theorem}[section]
\newtheorem{lemma}[theorem]{\indent\sc Lemma}

\theoremstyle{definition} 
\newtheorem{definition}[theorem]{\indent\sc Definition}
\newtheorem{remark}[theorem]{\indent\sc Remark}

\newcommand{\cb}{{\mathcal B}}

\newcommand{\vep}{\varepsilon}

\newcommand{\ssp}{\rm {supp}}
\newcommand{\setdef}{\stackrel {\rm {def}}{=}}
\newcommand{\beq}{\begin{equation}}
\newcommand{\eeq}{\end{equation}}
\newcommand{\R}{{\mathbb{R}}}
\newcommand{\Pro}{{\mathbb{P}}}
\newcommand{\E}{{\mathbb{E}}}
\newcommand{\T}{{\mathbb{T}}}

\newcommand{\Z}{{\mathbb{Z}}}
\newcommand{\N}{{\mathbb{N}}}

\def\build#1_#2^#3{\mathrel{\mathop{\kern 0pt#1}\limits_{#2}^{#3}}}
\def\tend#1#2{\build\hbox to 12mm{\rightarrowfill}_{#1\rightarrow #2}^{}}
\def\netendpas#1#2{\build\hbox to 12mm{$\not \longrightarrow$}_{#1 \rightarrow
		#2}^{}}

\def \converge#1#2#3{\build\hbox to 15mm {\rightarrowfill}_{#1\rightarrow #2}^
	{\hbox{\scriptsize #3}}}

\begin{document}

\title[On symmetric version of Saeki's Theorem]{On the symmetric version of Saeki's Theorem and flat densities} 

\author[e. H. El Abdalaoui ]{\lowercase{el} Houcein \lowercase{el} Abdalaoui$^*$} 



%
\thanks{ 
$^{*}$This  work was done while the author was in delegation at the CNRS, DR-19 Normandy. 
}
\address{
Department of Mathematics, LMRS, UMR 60 85,  \endgraf
University
of Rouen Normandy \endgraf
Avenue de l'Universit\'e, BP.12, 76801
Saint Etienne du Rouvray \endgraf
France 
}
\email{elhoucein.elabdalaoui@univ-rouen.fr}



\begin{abstract}
	It is shown that for any $\alpha \in ]\frac12,1[$ there exists a symmetric probability measure $\sigma$ on the torus such that
the Hausdorff dimension of its support is $\alpha$ and $\sigma*\sigma$ is absolutely continuous with flat continuous Radon-Nikodym derivative. Namely, we obtain a symmetric version of Saeki Theorem but the flat Radon-Nikodym derivative of $\sigma*\sigma$ can not be a Lipschitz function.
\end{abstract}
\maketitle
\section*{Introduction} 
Saeki's asserted that there exists a singular measure $\sigma$ with support of Lebesgue measure zero such that $\sigma*\sigma$  is absolutely continuous and in addition the Radon-Nikodym derivative of $\sigma*\sigma$ has a uniformly convergent Fourier
series \cite{Saeki}. It is an easy exercise to show that the measure $\sigma$ can not be symmetric. Nevertheless, using the ideas of K\"{o}rner's proof of Saeki Theorem \cite{Korner}, we shall obtain a symmetric version of Saeki theorem. In fact, we shall prove that there exist a symmetric singular measure for which the convolution is absolutely continuous with continuous Radon-Nikodym derivative. Our principal motivation is connected to the question raised in \cite{AL} on the existence of singular symmetric measure $\sigma$ with absolutely continuous convolution such that the Radon-Nikodym derivative is flat. We recall that the function $f$ is flat if $\big\|f-1\big\|_{\infty}<\vep,$ for some $\vep \in [0,1).$ The subject of this note is to establish that such measure exists. For that, we will essentially follow K\"{o}rner's proof of Saeki Theorem \cite{Korner}. \\

This note is organized as follows. In section 2,  we state our main result and the fundamental theorem which is used for its proof. In section 3, we present the ingredients which are needed in the proof of 
the fundamental theorem, and we conclude by presenting its proof.

\section{The main theorem and its proof}
	We start by stating our main theorem.
\begin{theorem}
\label{saeki} Given $\varepsilon>0$ we can find a symmetric measure $\sigma$ with support of Hausdorff dimension $\displaystyle \frac12$ and a continuous function $f~~~:~~~\T \longrightarrow \R$ such that $\big\|f-1\big\|_{\infty} < \varepsilon,$ and $\sigma *  \sigma = f d\lambda$, 
where $\lambda$ is the Lebesgue measure.
\end{theorem}
	The main ingredients of Theorem \ref{saeki} are contained in the K\"{o}rner's proof of Saeki Theorem \cite{Korner}. Following K\"{o}rner, we denote by $\mathcal{F}_s$ the space of non-empty closed symmetric subsets of $\T$ equipped with the Hausdorff distance $d_{H}$ define by
\[
\forall (E,F) \in \mathcal{F}_s^2,~~~d_{H}(E,F)=\sup_{e \in E}d(e,F)+\sup_{f \in F}d(E,f).
\]
It is an easy exercise to verify that $(\mathcal{F}_s,d_H)$ is a complete metric space, in fact, we have more $(\mathcal{F}_s,d_H)$ is compact (see \cite[Ch. IV.] {kura}). As in \cite{Korner}, we consider the metric space $(\mathcal{E}_s,d_{\mathcal{E}_s})$, where $\mathcal{E}_s$ is consisting of ordered pairs $(E,\mu)$ where $E \in \mathcal{F}_s$ and $\mu$ is a symmetric probability measure with $\ssp(\mu) \subset E$ and $\widehat{\mu}(r) \longrightarrow 0$ as $|r| \longrightarrow +\infty$ and $d_{\mathcal{E}_s}$is defined by
\[
\forall ((E,\mu),(F,\sigma)) \in \mathcal{E}_s^2,~~~{\mathcal{E}_s}((E,\mu),(F,\sigma))=d_{H}(E,F)+\sup_{r \in \Z}|\widehat{\mu}(r)-\widehat{\sigma}(r)|.
\]
Finally, we consider the metric space $(\mathcal{G}_s,d_{\mathcal{G}_s})$, where $\mathcal{G}_s$ 
consisting of those
$(E,\mu) \in \mathcal{E}_s$ such that $\mu * \mu =f_{\mu}d\lambda$ such that $f_{\mu}$ is continuous, $d_{\mathcal{G}_s}$ is given by, for all $((E,\mu),(F,\sigma)) \in \mathcal{G}_s$
\begin{align*}
d_{\mathcal{G}_s}((E,\mu),(F,\sigma))=d_{\mathcal{E}_s}((E,\mu),(F,\sigma))+
||f_{\mu}-f_{\sigma}||_{\infty}.
\end{align*}
We stress that the space  $\mathcal{G}_s$ is similar to the space $\mathcal{G}$ introduced in \cite[Lemma 18, (iii)]{Korner}.
\noindent For the proof of our main result, we need the following Lemma.
	\begin{lemma}\label{Comp} The metric spaces $(\mathcal{E}_s,d_{\mathcal{E}_s})$ and $(\mathcal{G}_s,d_{\mathcal{G}_s})$ are complete.
\end{lemma}
The proof of Lemma \ref{Comp} is similar to that of Lemma 18 in \cite{Korner} and it is leaved to the readers. We need also the following crucial lemma.
\begin{lemma}\label{Hn}
	Let $\alpha \in [\frac12,1)$ and $\mathcal{H}_{n}$ be the subset of consisting of those $(E,\mu) \in \mathcal{G}_s$ such that we can find a finite collection of intervals
	$\mathcal{I}$ symmetric (which means if $I \in {\mathcal{I}}$ the $-I$ is in ${\mathcal{I}}$) with
	\[
	E \subseteq \bigcup_{I \in \mathcal{I}}I {\rm {~~and~~}} \sum_{I \in \mathcal{I}}|I|^{\alpha+\frac1{n}}<\frac1{n}.
	\]
	Then $\mathcal{H}_{n}$ is an open dense set in $(\mathcal{G}_s,d_{\mathcal{G}_s})$.
\end{lemma}
\noindent  The proof of Lemma \ref{Hn} is similar to that of Lemma 24 in \cite{Korner}, since we can use Lemma 26 from \cite{Korner} without repeating its tedious proof, the only thing we need to show is an assertion corresponding to Lemma 25 in \cite{Korner}, for that we point out that the metric $d_{\psi}$ defined in \cite{Korner} satisfies
$d_{\psi} \geq d_{\mathcal{G}_s}$. Let us further emphasize that  Lemma \ref{Hn} plays the same role as Lemma 24 in \cite{Korner}.\\	
At this point, we state the fundamental result of this note.
\begin{theorem}\label{TKorner} Let $\alpha  \in [\frac12,1)$. The complement of the set
	$$
	\mathcal{H}_{\alpha}=\{(E,\mu) \in \mathcal{G}_s ~~:~~ E {\rm ~~has~~Hausdorff~~dimension~~} \alpha\}
	$$
	is of first category in $(\mathcal{G}_s,d_{\mathcal{G}_s})$.
\end{theorem}
\noindent Obviously, Theorem \ref{saeki} follows from Theorem \ref{TKorner}.
\section{The proof of Theorem \ref{TKorner}.}
The fundamental ingredient of the proof of Theorem \ref{TKorner} is based on the following lemmas.
The first one is standard in Probability and Martingale theory, and we provide its proof 
for the readers convenience. Our proof is different than that of its weaker version stated in  \cite{Korner}.
\begin{lemma}\label{Hoefding}
	Suppose that $0 < Np < 1$ and $ m \geq 2$. Then, if  $Y_1, Y_2,\cdots,$\\$ Y_N$
	are independent random variables with
	$$\Pro(Y_j = 1) = p,~~ \Pro(Y_j = 0) = 1-p,$$
	\noindent{}it follows that
	\[
	\Pro\left(\sum_{j=1}^{N}Y_j \geq m\right)
	\leq \frac{(Np)^m}{m!}.
	\]
\end{lemma}
\begin{proof}It is well known that the distribution of $S_N=\sum_{j=1}^{N}Y_j$ is a binomial distribution
	$\cb(N,p)$. Therefore
\begin{align}
	\Pro\left(S_N\geq m\right)&=\sum_{k=m}^{N}\cb(N,p)\big\{k\}\\
	&=\sum_{k=m}^{N}\binom{k}{N}p^k (1-p)^{N-k}\\
	&\leq  (Np)^m \sum_{k=m}^{N}\frac{(N-m)!}{k! (N-k)!}p^{k-m} (1-p)^{N-k} \label{Kor:eq1}\\
	&\leq \frac{(Np)^m}{m!} \sum_{k=m}^{N}\frac{(N-m)!}{(k-m)! (N-k)!}p^{k-m} (1-p)^{N-k} \label{Kor:eq2}
\end{align}
The inequality \eqref{Kor:eq1} is due to the fact that 
$$N! p^k= Np. ((N-1)p).\cdots.((N-m+1)p)\times (N-m)!,$$
and the last inequality \eqref{Kor:eq2} is a consequence of 
$$\binom{m}{k} \geq 1, \; \; \textrm{that~~is} \; \; k! \geq m! (k-m)!.$$
We thus get 
 \begin{align}
 \Pro\left(S_N\geq m\right)&=\sum_{k=m}^{N}\cb(N,p)\big\{k\}\\
 &\leq \frac{(Np)^m}{m!} \sum_{l=0}^{N-m}\frac{(N-m)!}{l! (N-m-l)!}p^{l} (1-p)^{N-m-l}\\
 &\leq  \frac{(Np)^m}{m!} \cb(N-m,p)\Big\{0,\cdots,N-m\Big\}= \frac{(Np)^m}{m!} 
 \end{align}
\end{proof}
\noindent Before stating the second one, we recall the following definition.
\begin{definition}
	A sequence $W_r$ is said to be a martingale with respect to a
	sequence $X_r$ of random variables if
	\begin{itemize}
		\item[(i)] $\E \left(\left|W_j\right|\right) <\infty.$
		\item[(ii)] $\E \left(W_{r+1}|X_0,\cdots,X_r\right)=W_r.$
	\end{itemize}
\end{definition}
For $n \geq 2$, we put 
	$$ \Gamma_n \setdef \left \{\frac{r}{n}~~:~~ r\in \left \{-1,\cdots,\frac{-1}{n},\frac{1}{n},\cdots,1\right \}\right\}.$$
\begin{lemma}\label{martingale}
	Let $\delta > 0$ and let $W_r$ be a martingale with respect to a sequence
	$X_r$ of random variables. Write $Y_{r+1} = W_{r+1}-W_r$. Suppose that
	$$\E\left(e^{\lambda Y_{r+1}}|X_0,X_1,\cdots,X_r\right) \leq e^{a_{r+1}\frac{\lambda^2}{2}}.$$
	for all $\lambda < \delta$ and some $a_{r+1} >0$. Suppose further that $A \geq \sum_{r=1}^{N}a_r$. Then,
	provided that $0 \leq x < A\delta$, we have
	$$
	\Pro\left(|W_N-W_0|\geq x \right) \leq \exp\Big(\frac{-x^2}{2A}\Big).
	$$
\end{lemma}
Lemma \ref{martingale} is known as Hoeffding-Azuma's inequality \cite[p.305]{RO}, \cite[p.73]{DP}, \cite[p. 237, E14.2]{Williams}. By applying Lemma \ref{Hoefding}, we get the following lemma. For its complete proof, we refer to \cite{Korner}.
\begin{lemma}\label{laborne}
	Let $\gamma \in ]0,1[$ and $\varepsilon>0$, we can find an integer $M=M(\gamma,\varepsilon)\geq 1$ such that the following property holds. Suppose
	$n \geq 2, n^{\gamma} \geq N$ and $X_1,X_2,\cdots,X_{N}$ are independent symmetric random variables each uniformly distributed on
	$$\Gamma_n \setdef \left \{\frac{r}{n}~~:~~ r\in \left \{-1,\cdots,\frac{-1}{n},\frac{1}{n},\cdots,1\right \}\right\}.$$
	Then, with probability at least $1-\frac{\varepsilon}{n}$,
	$$\sum_{j=1}^{N}\left(\delta_{X_j}\left ( \left \{ \frac{r}{n}\right\}\right)+\delta_{-X_j}\left ( \left \{ \frac{r}{n}\right\}\right)\right) < M.$$
	for all $r\in \left \{-1,\cdots,\frac{-1}{n},\frac{1}{n},\cdots,1\right \}.$
\end{lemma}
\noindent The key lemma is the following lemma.
\begin{lemma}\label{Kaumfman}
	Suppose $\varphi~~:~~ \N \longrightarrow \R$ is a sequence with $\varphi(n)\tend{n}{\infty}+\infty.$
	If $\gamma \in ]0,1[$ and $\varepsilon> 0$, there exists two integers $M\setdef M(\gamma)$ and $n_0(\varphi,\gamma,\varepsilon)\setdef n_0$
	with the following property. Suppose that $n >n_0$, $n$ is odd, $n^{\gamma} \geq N$ and $X_1,X_2,\cdots,X_N$ are independent symmetric random variables each uniformly distributed on $\Gamma_n$. Then, if we write $\sigma \setdef \frac{1}{2N}\sum_{j=1}^{N}\left (\delta_{X_j}+\delta_{-X_j}\right)$, we have
	$$\left|\sigma*\sigma\left(\frac{k}{n}\right)-\frac{1}{2n}\right| \leq \varepsilon \frac{\varphi(n)\sqrt{ln(n)}}{N\sqrt{n}},$$
	\noindent{}and
	$$\sigma\left(\frac{k}{n}\right) \leq \frac{M}{N}$$
	for all $k \in \left \{-1,\cdots,\frac{-1}{n},\frac{1}{n},\cdots,1\right \},$ with probability at least $\frac12$.
\end{lemma}
\noindent As a corollary of the lemme \ref{Kaumfman}, we have the following lemma.

\begin{lemma}\label{Kaumfman2}
	Suppose $\varphi~~:~~ \N \longrightarrow \R$ is a sequence with $\varphi(n)\tend{n}{\infty}+\infty.$
	If $\gamma \in ]0,1[$ and $\varepsilon> 0$, there exist two integer $M\setdef M(\gamma)$ and $n_0(\varphi,\gamma,\varepsilon)\setdef n_0$
	with the following property. Suppose that $n >n_0$, $n$ is odd, $n^{\gamma} \geq N$, we can find $N$ points
	$$x_j \in \Gamma_n,$$
	\noindent{} such that writing $\sigma \setdef \frac{1}{2N}\sum_{j=1}^{N}\left (\delta_{x_j}+\delta_{-x_j}\right)$, we have
	$$\left|\sigma*\sigma\left(\frac{k}{n}\right)-\frac{1}{2n}\right| \leq \varepsilon \frac{\varphi(n)\sqrt{ln(n)}}{N\sqrt{n}},$$
	\noindent{}and
	$$\sigma\left(\frac{k}{n}\right) \leq \frac{M}{N}$$
	for all $k \in \left \{-1,\cdots,\frac{-1}{n},\frac{1}{n},\cdots,1\right \}$.
\end{lemma}
\noindent Now, let us emphasize that we need only to give a sketch of the proof of the lemma \ref{Kaumfman}. We further point out that the outline of the proof follows the proof of Lemma 10 in \cite{Korner}, but here the definition of the martingale difference $Y_j$ is more complicated. 
\begin{proof} Let $M=M(\gamma,\frac14)$ be as in Lemma \ref{laborne}. Fix $\frac{r}{n} \in \Gamma_n$ and define $Y_1,Y_2,\cdots,Y_N$ as follows. If
	$\displaystyle \sum_{v=1}^{j-1}\left (\delta_{X_v}\left(\frac{u}{n}\right)+\delta_{-X_v}\left(\frac{u}{n}\right)\right)<M,$ for all $u$ with $1 \leq |u| \leq n,$ set
	\newpage
	\begin{eqnarray*}
		Y_j =&&-\frac{2j-1}{2n}+\frac14\left(\delta_{2X_j}\left(\left\{\frac{r}{n}\right\}\right)+\delta_{-2X_j}\left(\left\{\frac{r}{n}\right\}\right)\right)+\\
		&&\frac14\sum_{v=1}^{j-1}\left\{\delta_{X_v+X_j}\left(\left\{\frac{r}{n}\right\}\right)+\delta_{-(X_v+X_j)}\left(\left\{\frac{r}{n}\right\}\right)\right\}+\\
		&&\frac14\sum_{v=1}^{j-1}\left\{\delta_{X_v-X_j}\left(\left\{\frac{r}{n}\right\}\right)+\delta_{X_j-X_v}\left(\left\{\frac{r}{n}\right\}\right)\right\}.
	\end{eqnarray*}
	\noindent{}Otherwise $Y_j=0$. Put $W_0=0$ and $W_j=\displaystyle \sum_{v=1}^{j}Y_v.$ It follows that 
	If $(X_1, \cdots, X_{j-1} )$ are given, then, each of the measures $\delta_{X_v+X_j}, \delta_{-\big(X_v+X_j\big)}, \delta_{X_v-X_j}$ and $\delta_{X_j-X_v}$ is uniform on $\Gamma_n$. Since $n$ is odd, each of the measures $\delta_{2X_j}$ and $\delta_{-2X_j}$ is also uniform on
	$\Gamma_n$ .  Hence,
\begin{align*}
\E\Big(\delta_{\pm \big(X_v+X_j\big)}\big\{\frac{r}{n}\big\}\Big|_{X_1,\cdots,X_{j-1}}\Big)&=
\E\Big(\delta_{\pm \big(X_v-X_j\big)}\big\{\frac{r}{n}\big\}\Big|_{X_1,\cdots,X_{j-1}}\Big)\\
&=\E\Big(\delta_{\pm 2X_j}\big\{\frac{r}{n}\big\}\Big|_{X_1,\cdots,X_{j-1}}\Big)=\frac{1}{n}.
\end{align*}

We thus get $\E(Y_j|X_1,\cdots,X_{j-1})=0$, and thus the sequence $W_j$ is a martingale with respect to $X_1,\cdots,X_N$. Following K\"{o}rner proof (of Lemma 10 in \cite{Korner}), we further get that
	$$\E(e^{\lambda Y_j}|X_1,\cdots,X_N) \leq \exp\left\{\frac{N}{n}4(1+M^2)\lambda^2\right\}.  $$
	We can thus apply Lemma \ref{martingale} with
	$$A=8\frac{N^2}{n}(M^2+1) ~~~{\rm{and}}~~~x=\varepsilon\frac{N\phi(n)\sqrt{\ln(n)}}{\sqrt{n}},$$
	\noindent since $\phi(n) \tend{n}{+\infty}{+\infty}$ we can choose $n_0(\phi,\gamma,\varepsilon)=n_0.$ Therefore 
	$$
	\Pro\left\{|W_N| \geq \varepsilon\frac{N\phi(n)\sqrt{\ln(n)}}{\sqrt{n}}\right\} \leq \frac1{4n},
	$$
	\noindent{}for all $n \geq n_0.$
	To finish the proof, observe that
	$$ |W_N|=\left|\sum_{j=1}^{N}Y_j\right|=\left|\left(N\sigma*N\sigma\right) \left(\left\{\frac{r}{n}\right\}\right)-\frac{N^2}{2n}\right|.$$
	It follows with probability at least $1-\frac1{2n}$, that we have
	$$\sum_{j=1}^{N}\left(\delta_{X_j}+\delta_{-X_j}\right)\left(\left\{\frac{r}{n}\right\}\right)<M,$$
	for all $r$ with $1 \leq |r| \leq n$ and
	$$|\left(N\sigma*N\sigma\right) \left(\left\{\frac{r}{n}\right\}\right)-\frac{N^2}{2n}| < \varepsilon\frac{N\phi(n)\sqrt{\ln(n)}}{\sqrt{n}},$$
	hence
	$$|\left(\sigma*\sigma\right) \left(\left\{\frac{r}{n}\right\}\right)-\frac{1}{2n}|<
	\varepsilon\frac{\phi(n)\sqrt{\ln(n)}}{N\sqrt{n}}.$$
	and the proof of the lemma is complete.
\end{proof}
Lemmas 2.5 and 2.6 play the same roles as Lemmas 10 and 9 in \cite{Korner} respectively.
Now we rewrite the lemma \ref{Kaumfman2} in more usable form. More precisely, we will exhibit a function $g$ in order to prove the following lemma.
\begin{lemma}Suppose $\varphi~~:~~ \N \longrightarrow \R$ is a sequence with $\varphi(n)\tend{n}{\infty}+\infty.$
	If $\gamma \in ]0,1[$ and $\varepsilon> 0$, there exists a integer $M(\gamma)\setdef M$ and $n_0(\varphi,\gamma,\varepsilon)\setdef n_0$
	with the following property. Suppose that $n >n_0$, $n$ is odd, $n^{\gamma} \geq N$, we can find $N$ points
	$$x_j \in \Gamma_n,$$
	\noindent{} such that, writing
	$$g=\frac{n}{N}\sum_{1 \leq |j| \leq N}\mathbbm{1}_{[x_j-\frac1{4n},x_j+\frac1{4n}]},$$
	with $x_{(-j)}=-x_j$, we have $g*g$ continuous and
	\begin{enumerate}
		\item $||g*g-1||_{\infty} \leq 2\varepsilon\frac{\phi(n)\sqrt{n\ln(n)}}{N}.$
		\item $|g(t)|\leq \frac{2nM}{N}$ for all $t \in \T$.
	\end{enumerate}
\end{lemma}
\begin{proof}
	Observe that we have
	$$g=\sigma*2n\mathbbm{1}_{[-\frac1{2n},\frac1{2n}]}.$$
	\noindent{}It follows that
	$$g*g=\sigma*\sigma*2n\mathbbm{1}_{[-\frac1{4n},\frac1{4n}]}*2n\mathbbm{1}_{[-\frac1{4n},\frac1{4n}]}=\sigma*\sigma*2n\Delta_n,$$
	\noindent{}where $\Delta_n(x)=\max\left\{0,1-2n|x|\right\}$. By the way, we get
	$$g*g(\frac{r}{n})=2n\left(\sigma*\sigma\right) \left(\left\{\frac{r}{n}\right\}\right).$$
\end{proof}
\noindent From now the rest of the proof follows the path of K\"{o}rner's proof and this finishes the proof of Theorem \ref{TKorner} and the proof of the main result of this note is done.
\begin{remark}Lemma 25 from \cite{Korner} tells us that we can approximate uniformly a continuous function by a sequence of infinitely differentiable functions but not in any space $\Lambda_{\psi}$ where $\psi$ is a positive strictly increasing continuous function which satisfies
	$$ \psi(0)=0 ~~{~\rm {and}}~~\frac{\psi(t)}{t^{\beta}}\tend{t}{0^{+}}0,~~\beta < \alpha-\frac12,$$
	\noindent{}for some given (in advance) $\alpha \in ]\frac12,1[.$
We recall that $\Lambda_{\psi}$  is the set of continuous
functions $f : \mathbb{T} \longrightarrow \mathbb{C}$ with 
$$\sup_{t,h \in \mathbb{T}, h \neq 0} \frac{|f(t+h)-f(t)|}{\psi(|h|)}<\infty.$$
\end{remark}

\end{document}